\documentclass[11pt]{article}
\usepackage[usenames]{color}
\usepackage{latexsym,amsmath,amssymb,amscd,amsthm,mathrsfs}
\usepackage{bbold}
\usepackage{amsmath,amsfonts}
\usepackage{bm}
\usepackage{eucal}

\begin{document}

\newtheorem{theorem}{\bf Theorem}[section]
\newtheorem{definition}{\bf Definition}[section]
\newtheorem{proposition}{\bf Proposition}[section]
\newtheorem{example}{\bf Example}[section]
\newtheorem{corollary}{\bf Corollary}[section]
\newtheorem{remark}{\bf Remark}[section]
\newtheorem{lemma}{\bf Lemma}[section]
\newcommand{\halmos}{\quad\hfill\mbox{$\Box$}}
\def\l{\ell}
\def \Liminf{\mathop{\underline{\lim}}\limits}
\def \Limsup{\mathop{\overline{\lim}}\limits}
\def\1{\mbox{1\hspace{-.25em}I}}
\def\Ex{{\bf E}}
\def\Pb{{\bf P}}
\let\bb\mathbb
\def\UU{{\bb U}}
\def\XX{{\bb X}}
\def\WW{{\bb W}}
\def\VV{{\bb V}}
\def\BB{{\bb B}}
\def\DD{{\bb D}}
\def\CC{{\bb C}}
\def\NN{{\bb N}}
\def\RR{{\bb R}}
\def\KK{{\bb K}}
\def\F{\mathcal{F}}
\def\iunit{{\bf i}}
\def\sgn{{\rm sgn}}
\def\eps{\varepsilon}
\let\bm\boldsymbol 

\title{On nonlinear TAR processes and threshold estimation}

\author{P. \textsc{Chigansky},
{\small
The Hebrew University,
 Jerusalem,
Israel}\\
Yu. A. \textsc{Kutoyants},
{\small  Universit\'e du Maine, Le Mans, France}}
\date{ }
\maketitle


\begin{abstract}
We consider the problem of threshold estimation for autoregressive time series
with a ``space switching'' in the situation, when the regression is nonlinear
and the innovations have a smooth, possibly non Gaussian, probability density.
Assuming that the unknown threshold parameter is sampled from a continuous
positive prior density,  we  find the asymptotic distribution of the Bayes
estimator. As usually in the singular estimation problems, the sequence of Bayes estimators
is  asymptotically efficient, attaining the minimax risk lower bound.
\end{abstract}



\section{Introduction}

The simplest threshold autoregressive (TAR) process  is the  time series, generated by the recursion
$$
X_{j+1}=\rho _1X_j\,\1_{\left\{X_j<\theta \right\}}+
\rho_2X_j\,\1_{\left\{X_j\geq \theta \right\}}+\varepsilon _{j+1} ,\quad
\quad j=0,\ldots, n-1,
$$
where  $\varepsilon _j\sim{\cal N}\left(0,\sigma ^2\right)$ are i.i.d. random variables and $ \rho _1\not=\rho _2$ and
$\sigma ^2$ are known constants. The unknown {\it
  threshold} parameter $\theta \in \Theta =\left(\alpha ,\beta \right)$ is to be estimated from the data
$X^n=\left(X_0,X_1,\ldots,X_n\right)$.   This
model and some of its generalizations  has been extensively studied  during the last decades (see e.g.
\cite{Ch93}-\cite{DN11},\cite{Tong11} and the references therein). Particularly,
much attention focused on  the properties of the least squares (LS) estimator
$$
 \theta _n^*={\rm argmin}_{\theta \in \Theta }
\sum_{j=0}^{n-1}\left[X_{j+1}-\rho _1X_j\,\1_{\left\{X_j<\theta \right\}}
-\rho_2X_j\,\1_{\left\{X_j\geq \theta \right\}}\right]^2.
$$
Assuming that $|\rho_1|\vee |\rho_2|<1$ and thus that $(X_j)$ is geometric mixing with the unique invariant density $\varphi \left(x,\theta \right)$,
Chan \cite{Ch93} proved  consistency of $\theta _n^*$ with the rate $n$
(rather than $\sqrt{n}$ as in regular problems) and showed that the limit distribution is related to certain compound Poisson process (see \eqref{Z} below).
Note that if $\varepsilon_1\sim  {\cal
  N}\left(0,\sigma ^2\right)$, the LS estimator coincides with the maximum likelihood (ML) estimator.

This work continues the study of the Bayes  estimator for the TAR models, initiated  in
\cite{ChK11} and  developed further in \cite{ChK10} and  \cite{CKL1} (see also \cite{Kut10} for
the continuous time counterpart). We consider  the following
more general nonlinear TAR(1) model
\begin{equation}
\label{nG}
X_{j+1}=h \left(X_j\right)\,\1_{\left\{X_j< \theta
\right\}}+g\left(X_j\right)\,\1_{\left\{X_j\geq  \theta
\right\}} +\varepsilon _{j+1} , \qquad j=0,\ldots, n-1,
\end{equation}
where $h(x)$ and $g(x)$ are known functions, $(\varepsilon _j)$ are
i.i.d. random variables with a known density function $f\left(x\right)>0,x\in
\RR$ and the initial condition $X_0$ is independent of $(\varepsilon_j)$ and
has a probability density $f_0(x)$.

\medskip
Throughout we shall assume that the following conditions are in force
\begin{enumerate}
\renewcommand{\theenumi}{({\bf a\arabic{enumi}})}
\renewcommand*\labelenumi{\theenumi}
\item\label{c1}  The parameter $\theta \in \left(\alpha ,\beta
  \right)\equiv \Theta , -\infty<\alpha <\beta <\infty $ is sampled from the continuous positive prior
  density $p\left(\theta \right),\theta \in\Theta $.
\item\label{c2} The functions $h$  and $g$ are
  continuous and  satisfy
$$
\inf_{v\in\Theta }\left|\delta \left(v\right)\right|>0,\qquad  \delta
\left(v\right):=g\left(v\right)-h\left(v\right).
$$

\item\label{c3} The random variables $(\varepsilon _j)_{j\ge 1}$
are i.i.d. with a known continuous bounded density function $f\left(x\right)>0,x\in \RR$

\item\label{c4} The functions $h\left(x\right),g\left(x\right)$ and
  $f\left(x\right)$ are such that the time series, generated by \eqref{nG}, is
  geometric mixing with the unique positive bounded invariant density $\varphi
  \left(x,\theta\right)$, i.e. for any measurable function $|\psi(x)|\le 1$
$$
\Ex\bigg|\Ex \big(\psi(X_j)|\F_i\big)-\int_\RR \psi(x)\varphi(x,\theta)dx\bigg|\le C r^{|j-i|},\quad j>i
$$
with positive constants $C$ and $r<1$.

\item\label{c5} The function
$$
J(z):=\int_{-\infty }^{\infty
}\left|\ln\frac{f\left(y+z\right)}{f\left(y\right)}\right|
f\left(y\right)\,{\rm d}y, \quad \min_{\theta\in \Theta} \delta(\theta)\le z\le \max_{\theta\in \Theta} \delta(\theta)
$$
is bounded.

\end{enumerate}

\medskip

The likelihood function of the sample $X^n$ is given by
$$
 L\left(\theta ,X^n\right)=f_0\left(X_0\right)\prod_{j=0}^{n-1}
f\Big(X_{j+1}-h \left(X_j\right)\,\1_{\left\{X_j< \theta
  \right\}}-g\left(X_j\right)\,\1_{\left\{X_j\geq \theta \right\}} \Big),
$$
and the Bayes estimator $\tilde\theta _n$ with respect to the mean square risk is the
conditional expectation
$$
\tilde\theta _n=\Ex\left(\theta | X^n\right)=\frac{\int_{\Theta }^{}\theta \,p\left(\theta \right)\,
  L\left(\theta ,X^n\right)\;{\rm d}\theta}{\int_{\Theta }^{}p\left(\theta
  \right)\, L\left(\theta ,X^n\right)\;{\rm d}\theta}.
$$
Since the likelihood $L(\theta,X^n)$ is piecewise constant in $\theta$,
the estimate can be computed efficiently (see \cite{ChK11}).

The asymptotic properties of $(\tilde\theta_n)$ are formulated in terms of
the following compound Poisson process
\begin{equation}\label{Z}
{Z}\left(u\right)= \left\{\begin{array}{ll}
&\exp\left( \,\sum_{l=1}^{N_+\left(\;u\,\right)}\,\, \ln
\frac{f\left(\varepsilon_l^++\delta \left(\theta _0\right)
  \right)}{f\left(\varepsilon_l^+ \right)} \right),\qquad u\geq 0 ,
\\ &\exp\left( \,\sum_{l=1}^{N_-\left(-u\right)} \ln
\frac{f\left(\varepsilon_l^--\delta \left(\theta _0\right)
  \right)}{f\left(\varepsilon_l^- \right)} \right),\qquad u< 0.
             \end{array}
             \right.
\end{equation}
Here $\theta_0  $ is the true value of the parameter, $\epsilon_l^\pm $
are independent random variables with the density function $f\left(x\right)$,
$N_+\left(\cdot \right)$, $N_-\left(\cdot \right)$ are independent Poisson
processes with the same intensity  $\lambda =\varphi \left(\theta
_0,\theta _0\right)$ ($Z(u):=1$ on the sets $\{N_\pm(u)=0\}$).

Define the random variable
$$
\tilde{u}=\frac{\int_{R}^{}u\,Z\left(u\right)\,{\rm d}u}{\int_{R}^{}Z\left(u\right)\,{\rm d}u}.
$$

As shown in \cite{ChK10} (see \cite{IH} for the general theory), we have the
following lower bound on the mean square risk of an arbitrary sequence of
estimates $(\bar\theta _n)$:
$$
\Liminf_{\delta \rightarrow 0}\Liminf_{n\rightarrow \infty
}\sup_{\left|\theta -\theta _0\right|<\delta } n^2\Ex_\theta
\left(\bar\theta _n-\theta \right)^2\geq \Ex_{\theta _0}\tilde u^2,
$$
and the Bayes estimates $(\widetilde \theta_n)$ are {\em efficient},  attaining this lower bound asymptotically.
Our main result is the following
\begin{theorem}
\label{T1}
Under the conditions \ref{c1}-\ref{c5}, the sequence of estimates
$(\tilde \theta _n)$ is consistent,  the convergence in distribution
$$
n\left(\tilde\theta _n-\theta _0\right)\Longrightarrow \tilde u
$$
holds and the moments converge:
$$
\lim_{n\rightarrow \infty }n^{p}\Ex_{\theta _0}\left|\tilde\theta
_n-\theta _0\right|^p=\Ex_{\theta _0}\left|\tilde u\right|^p, \quad p>0.
$$

\end{theorem}

\begin{remark}
The assumption \ref{c4} is often easy to check, using the standard ergodic theory as e.g. in \cite{MT09}.
The assumption \ref{c5}  is satisfied for many common densities. For example,
for the Gaussian innovations $\varepsilon _j\sim N\left(0,\sigma ^2\right)$,
$$
J\left(z\right)\le \frac{z^2}{2\sigma ^2}+\frac{|z|}{\sigma }.
$$
In this case, the limit compound Poisson process $Z\left(u\right)$ has Gaussian
jumps:
$$
\ln
\frac{f\left(\varepsilon_1^\pm\pm \delta \left(\theta _0\right)
  \right)}{f\left(\varepsilon_1^\pm \right)}
 =-\frac{\delta \left(\theta_0 \right)^2}{2\sigma ^2} \mp
\frac{\delta \left(\theta_0 \right)}{\sigma ^2}\;\varepsilon^\pm_1\sim{\cal
  N}\left( -\frac{\delta ^2\left(\theta _0\right)}{2\sigma ^2},\frac{\delta
  ^2\left(\theta _0\right)}{\sigma ^2} \right).
$$

Similarly the assumption  \ref{c5} is checked for
the Laplace density  $f\left(y\right)=\left(2\sigma
\right)^{-1}e^{-\frac{\left|y\right|}{\sigma }}$
and the  limit process has jumps of the form
$$
\ln
\frac{f\left(\varepsilon_1^\pm\pm\delta \left(\theta _0\right)
  \right)}{f\left(\varepsilon_1^\pm \right)}
 =\frac 1 {\sigma^2} \Big(\big|\eps^{\pm}_1\big|-\big|\eps^{\pm}_1\pm\delta(\theta_0)\big|\Big).
$$

\end{remark}

\section{The Proof}
We shall verify the conditions of the Theorem 1.10.2 in  \cite{IH},
where the properties of the Bayes estimators, announced in Theorem \ref{T1},  are derived from the convergence
of the normalized likelihood ratios
$$
Z_n\left(u\right)=\frac{ L\left(\theta_0+u/n,X^n \right) }{
  L\left(\theta_0,X^n\right)},\qquad u\in \UU_n=\left[n\left(\alpha
  -\theta _0\right),n\left(\beta -\theta _0\right)\right]
$$
to the limit process $Z\left(u\right),u\in \RR$  and the two inequalities
\eqref{e1} and \eqref{e2},  presented  below. The change of
 variables $\theta =\theta _0+u/n$ gives
$$
\tilde\theta _n=\frac{\int_{\UU_n}^{}\left(\theta
  _0+\frac{u}{n}\right)\, p\left(\theta _0+\frac{u}{n}\right)
\frac{L\left(\theta _0+\frac{u}{n},X^n \right)}{L\left(\theta
  _0,X^n\right)} {\rm d}u} {\int_{\UU_n}^{} p\left(\theta _0+\frac{u}{n}\right)
\frac{L\left(\theta _0+\frac{u}{n},X^n \right)}{L\left(\theta
  _0,X^n\right)} {\rm d}u}
=\theta
_0+\frac{1}{n}\frac{\int_{\UU_n}u\left[p\left(\theta
    _0\right)+o\left(1\right)\right] Z_n\left(u\right){\rm d}u}{\int_{\UU_n}\left[p\left(\theta
    _0\right)+o\left(1\right)\right] Z_n\left(u\right){\rm d}u}.
$$
Then, informally, we have
\begin{align*}
\tilde u_n=n\left(\tilde\theta _n-\theta
_0\right)=\frac{\int_{\UU_n}u\;Z_n\left(u\right){\rm
    d}u}{\int_{\UU_n}Z_n\left(u\right){\rm d}u}+o\left(1\right)\Longrightarrow
\frac{\int_{\RR}u\;Z\left(u\right){\rm
    d}u}{\int_{\RR}Z\left(u\right){\rm d}u}\equiv \tilde u.
\end{align*}
Theorem 1.10.2 in \cite{IH} validates  this convergence along  with the
convergence of moments. Similar program has been  realized  in the preceding works
\cite{ChK11}, \cite{ChK10} and \cite{CKL1}.

\begin{remark}
To avoid inessential technicalities, we shall assume that $(X_j)$ is
stationary, i.e.  $X_0\sim \varphi(\cdot,\theta_0)$. Due to the mixing
property \ref{c4}, all the results below can be derived without stationarity
assumption, along the same lines with minor adjustments (see \cite{CKL1} for
details).
\end{remark}

\begin{remark}
Below, $C$, $C'$, $c$, $C_p$, etc. denote constants, whose values are not important and
may change from line to line.  We shall denote by $\Pb_{\theta}$ and
$\Ex_{\theta}$ the probability and the expectation, corresponding to the
particular value of the unknown parameter $\theta\in \Theta$ and set
$\F_j:=\sigma\{\eps_i, i\le j\}$. The standard $O(\cdot)$ and $o(\cdot)$
notations will be used and we set $\sum_{i=k}^m(...)=0$ and $\prod_{i=k}^m
(...)=1$ for $k>m$.
\end{remark}
\subsection{Convergence of f.d.f.}

We shall prove the convergence of the finite dimensional distributions:
\begin{equation}
\label{wc1}
\big(\ln Z_n(u_1),...,\ln Z_n(u_d)\big)
\Longrightarrow
\big(\ln Z(u_1),...,\ln Z(u_d)\big), \quad u\in \RR^d,
\end{equation}
following \cite{CKL1}. We shall restrict the consideration to
$0=u_0<u_1<...<u_d$, leaving out the similar complementary case.  To this end,
note that the declared limit process $\ln Z(u)$ has independent increments and
\begin{multline*}
\Ex_{\theta_0} \exp\left(\sum_{j=1}^d \iunit\lambda_j \Big(\ln Z(u_j)-\ln Z(u_{j-1})\Big) \right) =\\
\exp \left(
\sum_{j=1}^d (u_j-u_{j-1})\varphi(\theta_0,\theta_0) \Big(\Psi(\lambda_j)-1\Big)
\right)=: e^{H(\lambda)}, \quad \lambda \in \RR^d,
\end{multline*}
where  (recall that $\delta:=g-h$)
$$
\Psi(\lambda_j) := \Ex_{\theta_0} \exp\left(\iunit \lambda_j \ln
\frac{f\big(\varepsilon_1+\delta \left(\theta _0\right)
  \big)}{f\big(\varepsilon_1 \big)}\right).
$$
Since $\ln Z(0)=0$ a.s., \eqref{wc1} follows from the convergence of characteristic functions
of the increments
$$
\lim_n \Ex_{\theta_0} \exp\left(\sum_{j=1}^d \iunit\lambda_j \Big(\ln Z_n(u_j)-\ln Z_n(u_{j-1})\Big) \right)
=
e^{H(\lambda)}, \quad \lambda \in \RR^d.
$$
Let $m(x,\theta):= h(x)\1_{\{x<\theta\}}+g(x)\1_{\{x\ge \theta\}}$ and note that
$$
m(x,\theta_0+u_{j-1}/n)-m(x,\theta_0+u_j/n) = \delta(x)\1_{\{x\in \DD^n_j\}},
$$
where $\DD^n_j:= [\theta_0+u_{j-1}/n, \theta_0+u_j/n)$.
Let $\BB^n_{j-1}:= [\theta_0,\theta_0+u_{j-1}/n)$,
then
\begin{equation}\label{eln}
\begin{aligned}
\ln Z_n(u_j)-&\ln Z_n(u_{j-1}) = \sum_{k=0}^{n-1}
\ln \frac
{
f\Big(X_{k+1}-m(X_k,\theta_0+u_j/n)\Big)
}
{
f\Big(X_{k+1}-m(X_k,\theta_0+u_{j-1}/n)\Big)
}= \\
&
\sum_{k=0}^{n-1}
\ln \frac
{
f\Big(\eps_{k+1}+m(X_k, \theta_0)-m(X_k,\theta_0+u_j/n)\Big)
}
{
f\Big(\eps_{k+1}+m(X_k, \theta_0)-m(X_k,\theta_0+u_{j-1}/n)\Big)
}=\\
&
\sum_{k=0}^{n-1}
\ln \frac
{
f\Big(\eps_{k+1}+\delta(X_k)\1_{\{X_k\in \BB^n_{j-1}\}}+\delta(X_k)\1_{\{X_k\in \DD^n_j\}})\Big)
}
{
f\Big(\eps_{k+1}+\delta(X_k)\1_{\{X_k\in \BB^n_{j-1}\}}\Big)
}=\\
&
\sum_{k=0}^{n-1}
\1_{\{X_k\in \DD^n_j\}}\ln \frac
{
f\Big(\eps_{k+1}+\delta(X_k)\1_{\{X_k\in \BB^n_{j-1}\}}+\delta(X_k)\Big)
}
{
f\Big(\eps_{k+1}+\delta(X_k)\1_{\{X_k\in \BB^n_{j-1}\}})\Big)
}\stackrel{\dagger}{=}\\
&
\sum_{k=0}^{n-1}
\1_{\{X_k\in \DD^n_j\}}\ln \frac
{
f\big(\eps_{k+1}+\delta(X_k)\big)
}
{
f\big(\eps_{k+1}\big)
}
=: \sum_{k=0}^{n-1} s^j_{k}
\end{aligned}
\end{equation}
where the equality $\dagger$ holds $\Pb_{\theta_0}$-a.s., since
$\Pb_{\theta_0}\Big(X_{k-1}\in \BB^n_{j-1}\cap \DD^n_j\Big)=0$.  Further,
define
$$
S_n:= \sum_{j=1}^d \lambda_j \Big(\ln Z_n(u_j)-\ln Z_n(u_{j-1})\Big)=\sum_{j=1}^d \lambda_j \sum_{k=0}^{n-1} s^j_{k}.
$$ We shall partition $n$ terms of this sum into $n^{1/2}$ consecutive blocks
of size $n^{1/2}$ and discard from each block its $n^{1/4}$ first entries. As
we shall see, this does not alter the asymptotic distribution of $S_n$, but
makes the blocks almost independent. Since in each block, the single event
$\{X_k\in D^n_j\}$ occurs with probability of order $n^{1/2}$, the Poisson
behavior emerges.  To implement these heuristics, define
$$
S_{m,n} := \sum_{j=1}^d \lambda_j \sum_{k = (m-1) n^{1/2}+n^{1/4}}^{m n^{1/2}} s_k^j, \quad m= 1,...,n^{1/2},
$$
and set
$
\widetilde S_n:= \sum_{m=1}^{n^{1/2}}S_{m,n}
$ (this is the sum, in which the $n^{1/4}$ entries of each block have been discarded).
By the triangle inequality
\begin{multline}\label{tri}
\Big|\Ex_{\theta_0} e^{\iunit S_n}- e^{H(\lambda)}\Big|\le
\Big|\Ex_{\theta_0} e^{\iunit S_n}- \Ex_{\theta_0} e^{\iunit \widetilde S_n}\Big| + \\
\bigg|\Ex_{\theta_0} e^{\iunit \widetilde S_n} -\Big(\Ex_{\theta_0} e^{\iunit S_{1,n}}\Big)^{n^{1/2}}\bigg|
+
\bigg|\Big(\Ex_{\theta_0} e^{\iunit S_{1,n}}\Big)^{n^{1/2}}- e^{H(\lambda)}\bigg|.
\end{multline}
We shall show that all the terms on the right hand side vanish as $n\to\infty$. By stationarity and the
assumption  \ref{c5},
\begin{align*}
&\Big|\Ex_{\theta_0} e^{\iunit S_n}-\Ex_{\theta_0} e^{\iunit\widetilde
    S_n}\Big|\le \Ex_{\theta_0} \Big| e^{\iunit S_n}-e^{\iunit\widetilde
    S_n}\Big| \le
\Ex_{\theta_0} \Big|  S_n-\widetilde   S_n\Big|\le\\
&n^{3/4}  \max_j |\lambda_j|\Ex_{\theta_0} \1_{\{X_0\in \DD^n_j\}}
\left|
 \ln \frac
{
f\big(\eps_{1}+\delta(X_0)\big)
}
{
f\big(\eps_{1}\big)
}
\right| = \\
&
n^{3/4} \max_j |\lambda_j|\Ex_{\theta_0} \1_{\{X_0\in \DD^n_j\}} J\big(\delta(X_0)\big)\le \\
&n^{3/4} \max_j |\lambda_j| \frac{u_j-u_{j-1}}{n}\sup_{x\in \RR}\varphi(x,\theta_0)\sup_{\theta\in \Theta}J(\delta(\theta))
\xrightarrow{n\to\infty}0,
\end{align*}
i.e. the first term in \eqref{tri} converges to zero.

Further, note that by the Markov property of $(X_j)$ and  \ref{c4}
$$
\bigg|\Ex_{\theta_0} \Big(e^{\iunit  S_{\ell,n}}\big|\F_{(\ell-1)n^{1/2}}\Big)-\Ex_{\theta_0} e^{\iunit S_{1,n}}\bigg| \le C r^{n^{1/4}}, \quad
\ell=1,...,n^{1/2}
$$
and hence
\begin{align*}
&\bigg|\Ex_{\theta_0} e^{\iunit \widetilde S_n} -\Big(\Ex_{\theta_0} e^{\iunit S_{1,n}}\Big)^{n^{1/2}}\bigg| =
\bigg|\Ex_{\theta_0}\prod_{m=1}^{n^{1/2}} e^{\iunit S_{m,n}} -\Big(\Ex_{\theta_0} e^{\iunit S_{1,n}}\Big)^{n^{1/2}}\bigg|=\\
&\bigg|\sum_{\ell=1}^{n^{1/2}}
\bigg(
\Ex_{\theta_0}\prod_{m=1}^{\ell} e^{\iunit S_{m,n}}
\Big( \Ex_{\theta_0} e^{\iunit S_{1,n}} \Big)^{n^{1/2}-\ell}
-
\Ex_{\theta_0}\prod_{m=1}^{\ell-1} e^{\iunit S_{m,n}}
\Big( \Ex_{\theta_0} e^{\iunit S_{1,n}} \Big)^{n^{1/2}-\ell+1}
\bigg)
\bigg|=\\
&
\bigg|\sum_{\ell=1}^{n^{1/2}}
\bigg(
\Ex_{\theta_0}\prod_{m=1}^{\ell-1} e^{\iunit S_{m,n}}
\Big(
e^{\iunit S_{\ell,n}}-\Ex_{\theta_0} e^{\iunit S_{1,n}}
\Big)
\Big( \Ex_{\theta_0} e^{\iunit S_{1,n}} \Big)^{n^{1/2}-\ell}
\bigg)
\bigg|=\\
&
\bigg|\sum_{\ell=1}^{n^{1/2}}
\bigg(
\Ex_{\theta_0}\prod_{m=1}^{\ell-1} e^{\iunit S_{m,n}}
\Big(
\Ex_{\theta_0} \Big(e^{\iunit S_{\ell,n}}\big|\F_{(\ell-1)n^{1/2}}\Big)-\Ex_{\theta_0} e^{\iunit S_{1,n}}
\Big)
\Big( \Ex_{\theta_0} e^{\iunit S_{1,n}} \Big)^{n^{1/2}-\ell}
\bigg)
\bigg|\le \\
&
\sum_{\ell=1}^{n^{1/2}}\Ex_{\theta_0}
\bigg|
\Ex_{\theta_0} \Big(e^{\iunit S_{\ell,n}}\big|\F_{(\ell-1)n^{1/2}}\Big)-\Ex_{\theta_0} e^{\iunit S_{1,n}}
\bigg|\le C n^{1/2} r^{n^{1/4}}\xrightarrow{n\to\infty}0.
\end{align*}
It is left to show that the last term in \eqref{tri} converges to zero. Let $\DD^n = \bigcup_{j=1}^d \DD^n_j$
and introduce the following events
\begin{align*}
& A_0 := \bigcap_{\ell\le n^{1/2}}\{X_\ell\not\in \DD^n\},\quad A_1 :=
  \bigcup_{j=1}^d\bigcup_{\ell=0}^{n^{1/2}} A_{\ell,j},\qquad A_{2+} :=
  \Big(A_0\cup A_1\Big)^c \\
 & A_{k,j} := \{X_k\in \DD^n_j\}\cap
  \bigcap_{\ell\le n^{1/2}, \ell\ne k} \{X_\ell\not \in \DD^n\} .
\end{align*}
In words, $A_0$ is the event, on which none of the first $n^{1/2}$ samples
falls in any of $\DD^n_j$'s, $A_1$ is the event of having exactly single
sample visiting one of $\DD^n_j$'s, etc. On the event $A_{k,j}$,
$$
S_{1,n} = \sum_{i=1}^d \lambda_i \sum_{\ell =n^{1/4}}^{ n^{1/2}} \1_{\{X_\ell\in \DD^n_i\}}\ln \frac
{
f\big(\eps_{\ell+1}+\delta(X_\ell)\big)
}
{
f\big(\eps_{\ell+1}\big)
}=
 \lambda_j  \ln \frac
{
f\big(\eps_{k+1}+\delta(X_k)\big)
}
{
f\big(\eps_{k+1}\big)
}
$$
and, since $\{X_k\in \DD^n_j\}=A_{k,j}\biguplus\Big(\{X_k\in \DD^n_j\}\cap \bigcup_{\ell\ne k} \{X_\ell \in \DD^n\}\Big)$,
\begin{align}\label{plugme}
&\Ex_{\theta_0} e^{\iunit S_{1,n}}\1_{A_1} = \sum_{j=1}^d \sum_{k=0}^{n^{1/2}} \Ex_{\theta_0} e^{\iunit S_{1,n}}\1_{A_{k,j}}=
\sum_{j=1}^d \sum_{k=0}^{n^{1/4}-1} \Pb_{\theta_0}(A_{k,j})+\\
\nonumber
&\sum_{j=1}^d \sum_{k=n^{1/4}}^{n^{1/2}} \Ex_{\theta_0} \exp\bigg(\iunit \lambda_j  \ln \frac
{
f\big(\eps_{k+1}+\delta(X_k)\big)
}
{
f\big(\eps_{k+1}\big)
}\bigg)\Big(\1_{\{X_k\in \DD^n_j\}}-\1_{\{X_k\in \DD^n_j\}\cap \bigcup_{\ell\ne k} \{X_\ell \in \DD^n\}}\Big).
\end{align}
By continuity  of $\varphi(x,\theta_0)$ and $\delta(x)$,
$$
\Pb_{\theta_0}(A_{k,j})\le \Pb_{\theta_0}(X_k\in \DD^n_j)=\frac{u_j-u_{j-1}}{n} \varphi(\theta_0,\theta_0) + o(n^{-1}),
$$
and
\begin{align*}
&\Ex_{\theta_0} \exp\bigg(\iunit \lambda_j  \ln \frac
{
f\big(\eps_{k+1}+\delta(X_k)\big)
}
{
f\big(\eps_{k+1}\big)
}
\bigg)\1_{\{X_k\in \DD^n_j\}} =\\
&\Ex_{\theta_0} \exp\bigg(\iunit \lambda_j  \ln \frac
{
f\big(\eps_{k+1}+\delta(\theta_0)\big)
}
{
f\big(\eps_{k+1}\big)
}
\bigg)\1_{\{X_k\in \DD^n_j\}} + o(n^{-1})=\\
&\Psi(\lambda_j)\frac{u_j-u_{j-1}}{n}\varphi(\theta_0,\theta_0) + o(n^{-1}).
\end{align*}
Further, by the Markov property, for $k<\ell$
\begin{align*}
& \Pb_{\theta_0}\big(X_k\in \DD^n_j, X_\ell \in \DD^n\big) = \Ex_{\theta_0}\1_{X_k\in \DD^n_j}\Pb_{\theta_0}\big(X_\ell \in \DD^n|\F_{\ell-1}\big) =\\
&
\Ex_{\theta_0}\1_{X_k\in \DD^n_j}\int_{\DD^n}
f\Big(x-h(X_{\ell-1})\1_{\{X_{\ell-1}<\theta_0\}}-g(X_{\ell-1})\1_{\{X_{\ell-1}\ge\theta_0\}}\Big)dx \le \\
& C_1 n^{-1} \Pb_{\theta_0}\big(X_k\in \DD^n_j\big) \le C_2n^{-2},
\end{align*}
where the inequalities hold, since the density $f(x)$ and therefore the invariant density  $\varphi(x,\theta_0)$, $x\in \RR$ are bounded.
Similar bound holds for $k>\ell$ and
it follows that
$$
\Pb_{\theta_0}\bigg(\{X_k\in \DD^n_j\}\cap \bigcup_{\ell\ne k, \ell\le n^{1/2}} \{X_\ell \in \DD^n\}\bigg) \le
\sum_{\ell\le n^{1/2}, \ell\ne k}\Pb_{\theta_0}\big(X_k\in \DD^n_j, X_\ell \in \DD^n\big) \le  C_3 n^{-3/2}.
$$
Plugging these estimates into \eqref{plugme}, we get
$$
\Ex_{\theta_0} e^{\iunit S_{1,n}}\1_{A_1}=n^{-1/2}\sum_{j=1}^ d\Psi(\lambda_j)\big(u_j-u_{j-1}\big)\varphi(\theta_0,\theta_0)+o(n^{-1/2}).
$$
If we set all $\lambda_j$'s to zeros, we also obtain
\begin{equation}
\label{A1}
\Pb_{\theta_0}(A_1)=n^{-1/2}\sum_{j=1}^ d\big(u_j-u_{j-1}\big)\varphi(\theta_0,\theta_0)+o(n^{-1/2}).
\end{equation}
Further,
\begin{multline*}
\Pb_{\theta_0}(A_0)=1- \Pb_{\theta_0}\left(\bigcup_{\ell\le n^{1/2}}\{X_\ell\in \DD^n\}\right) \ge 1-
\sum_{\ell\le n^{1/2}}\Pb_{\theta_0}\big(X_\ell\in \DD^n\big)= \\
 1-
\sum_{\ell\le n^{1/2}}\sum_{j=1}^d\Pb_{\theta_0}\big(X_\ell\in \DD^n_j\big)=1-
n^{-1/2}\sum_{j=1}^ d\big(u_j-u_{j-1}\big)\varphi(\theta_0,\theta_0)+o(n^{-1/2}).
\end{multline*}
On the other hand, $\Pb_{\theta_0}(A_0)\le 1-\Pb_{\theta_0}(A_1)$ and in view of \eqref{A1}, it follows that
\begin{equation}
\label{A0}
\Pb_{\theta_0}(A_0)=1-
n^{-1/2}\sum_{j=1}^ d\big(u_j-u_{j-1}\big)\varphi(\theta_0,\theta_0)+o(n^{-1/2}).
\end{equation}
Finally, using \eqref{A1} and \eqref{A0}, we also have
$$
\Pb_{\theta_0}(A_{2+}) = 1-\Pb_{\theta_0}(A_0)-\Pb_{\theta_0}(A_1) = o(n^{-1/2}).
$$
Assembling all parts together, we obtain the asymptotic
\begin{multline*}
\Ex_{\theta_0} e^{\iunit S_{1,n}} = \Pb_{\theta_0}(A_0)+\Ex_{\theta_0} e^{\iunit S_{1,n}}\1_{A_1} + \Ex_{\theta_0}  e^{\iunit S_{1,n}}\1_{A_{2+}} =\\
1+n^{-1/2}\bigg(\sum_{j=1}^ d  \Big(\Psi(\lambda_j)-1\Big)\big(u_j-u_{j-1}\big)\varphi(\theta_0,\theta_0)
\bigg) +o(n^{-1/2}),
\end{multline*}
and, in turn,
$$
\lim_n \bigg|\Big(\Ex_{\theta_0} e^{\iunit S_{1,n}}\Big)^{n^{1/2}}- e^{H(\lambda)}\bigg|=0.
$$
The claim now follows from \eqref{tri}.

\subsection{Equicontinuity}

The next step is to show that for some $C>0$
\begin{equation}
\label{e1}
\Ex_{\theta _0}\left(Z_n^{1/2}\left(u_2\right)-Z_n^{1/2}\left(u_1\right)\right)^2\leq C\,\left|u_2-u_1\right|.
\end{equation}
As in Lemma 2.4 in \cite{CKL1},  for e.g.  $u_2>u_1>0$, \eqref{eln} gives
\begin{align*}
&\Ex_{\theta_0}\left(Z_n^{1/2}\left(u_2\right)-Z_n^{1/2}\left(u_1\right)\right)^2\leq
\Ex_{\theta _0+{u_1}/n} \ln
\frac{Z_n\left(u_1\right)}{Z_n\left(u_2\right)}\\
& \le \Ex_{\theta _0+{u_1}/n}
\sum_{j=0}^{n-1}\left|\ln\frac{f \left(\varepsilon _{j+1} \right)
}{f\left(\delta \left(X_j\right)+\varepsilon _{j+1}    \right)}\right| \,\1_{\left\{\theta
    _0 +\frac{u_1}{n}\leq X_j < \theta_0+\frac{u_2}{n} \right\}}\\
& =\Ex_{\theta _0+{u_1}/n}
\sum_{j=0}^{n-1}\Ex_{\theta _0+{u_1}/n} \left(\left|\ln\frac{f
    \left(\varepsilon _{j+1} \right)
}{f\left(\delta \left(X_j\right)+\varepsilon _{j+1}    \right)}\right|\Big|{\cal F}_j \right) \,\1_{\left\{\theta
    _0 +\frac{u_1}{n}\leq X_j < \theta_0+\frac{u_2}{n} \right\}}\\
& =\Ex_{\theta _0+{u_1}/n}
\sum_{j=0}^{n-1}J\left(\delta \left(X_j\right)\right) \,\1_{\left\{\theta
    _0 +\frac{u_1}{n}\leq X_j < \theta_0+\frac{u_2}{n} \right\}}\\
& =n\int_{\theta
    _0 +\frac{u_1}{n}}^{\theta_0+\frac{u_2}{n} }J\left(\delta
\left(x\right)\right) \,\varphi \left(x,\theta_0+u_1/n\right)\,{\rm d}x \leq C\,\left|u_2-u_1\right|,
\end{align*}
as required.

\subsection{Large deviations estimate}

Finally we shall prove that for any $p>0$ there exists a constant $C_p>0$
such that
\begin{equation}
\label{e2}
\Ex_{\theta _0}Z_n^{1/2}\left(u\right)\leq \frac{C_p}{\left|u\right|^p}.
\end{equation}
We shall only sketch the proof,  as most of the arguments can be directly adopted from the proof of
Lemma 2.2 in \cite{ChK11} or Lemma 2.5, \cite{CKL1}.
Note that for any $c>0$,
\begin{multline}\label{une}
\Ex_{\theta _0}Z_n^{1/2}(u) = \Ex_{\theta _0}Z_n^{1/2}(u)\1_{\{Z_n^{1/2}(u)>e^{-c|u|}\}}
+
\Ex_{\theta _0}Z_n^{1/2}(u)\1_{\{Z_n^{1/2}(u)\le e^{-c|u|}\}}\le \\
\big(\Ex_{\theta _0}Z_n(u)\big)^{1/2}\Pb_{\theta_0}^{1/2}\big(Z_n^{1/2}(u)>e^{-c|u|}\big)
+
e^{-c|u|} = \Pb_{\theta _0}^{1/2}\big(\ln
Z_n^{1/2}\left(u\right)>-c\left|u\right|\big) +e^{-c\left|u\right|}
\end{multline}
and hence it suffices to show that for some $c>0$,
$$
 \Pb_{\theta _0}\big(\ln
Z_n^{1/2}\left(u\right)>-c\left|u\right|\big)\le \frac{C_p}{\left|u\right|^p}, \quad p>0.
$$
For $u>0$ (and similarly for $u<0$),
\begin{multline}\label{doux}
\Pb_{\theta _0}\left\{\ln
  Z_n^{1/2}\left(u\right)>-c\left|u\right|\right\}\\
 =\Pb_{\theta
    _0}\left\{\sum_{j=0}^{n-1}\ln\left[\frac{f \left(\delta
      \left(X_j\right)+\varepsilon _{j+1} \right) }{f\left(\varepsilon _{j+1}
      \right)}\right] \,\1_{\left\{\theta _0 \leq X_j <
    \theta_0+u/n \right\}}>-2cu\right\}.
\end{multline}
Let
$
\ell\left(x,y\right):=\ln\left[\frac{f \left(\delta
      \left(x\right)+y \right) }{f\left(y
      \right)}\right]
$
and introduce the notations
\begin{align*}
 G\left(\delta \right)&=-\ln H\left(\delta\right),\\
S_n^{\left(1\right)}&=\sum_{j=0}^{n-1}\ell\left(X_j,\varepsilon
 _{j+1}\right) \,\1_{\left\{X_j\in\BB^n   \right\}},\quad
 S_n^{\left(2\right)}=\sum_{j=0}^{n-1}G\left(\delta\left(X_j\right)
 \right) \,\1_{\left\{X_j\in\BB^n \right\}},
\end{align*}
where $\BB^n=[\theta_0, \theta_0+u/n]$ and
$$
H\left(\delta\right):=\int_{-\infty }^{\infty }
 \left(\frac{f\left(\delta
   +y\right)}{f\left(y\right)}\right)^{1/2}f\left(y\right)\,{\rm
   d}y
$$
is the Hellinger integral of order $1/2$.
By the Jensen inequality for all $\delta\ne 0$, $H(\delta)<1$
and hence $G\left(\delta \right)>0$.

Further, we have the following identity
\begin{align}
\label{eg}
\Ex_{\theta _0}e^{\frac 1 2 S_n^{\left(1\right)}+S_n^{\left(2\right)}}=1.
\end{align}
Indeed
\begin{align*}
\Ex_{\theta
  _0}e^{\frac 1 2 S_n^{\left(1\right)}+S_n^{\left(2\right)}}=\Ex_{\theta
  _0}e^{\frac 1 2 S_{n-1}^{\left(1\right)}+S_{n-1}^{\left(2\right)}}\Ex_{\theta
  _0}\left(\left. e^{\left[\frac 1 2 \ell\left(X_{n-1},\varepsilon
 _{n}\right)+G\left(\delta
    \left(X_{n-1}\right)\right)\right]\1_{\left\{X_{n-1}\in\BB^n
    \right\}}}\right|{\cal F}_{n-1}\right)
\end{align*}
and
\begin{align*}
&\Ex_{\theta _0}\left(\left. e^{ \frac 1 2 \ell\left(X_{n-1},\varepsilon
    _{n}\right)\1_{\left\{X_{n-1}\in\BB^n \right\}}}\right|{\cal
    F}_{n-1}\right)=\\
&\qquad \Ex_{\theta _0}\left(\left. e^{\frac 1 2
    \ell\left(X_{n-1},\varepsilon _{n}\right)}\right|{\cal F}_{n-1}\right)
  \1_{\left\{X_{n-1}\in\BB^n
    \right\}}+\1_{\left\{X_{n-1}\not\in\BB^n\right\}} =\\
&\qquad
 \int_{-\infty }^{\infty }\left(\frac{f\left(\delta
    \left(X_{n-1}\right)+y\right)}{f\left(y\right)}\right)^{1/2}f\left(y\right){\rm
    d}y \;\1_{\left\{X_{n-1}\in\BB^n\right\}}+\1_{\left\{X_{n-1}\not\in\BB^n \right\}}=\\
&\qquad
  \exp\Big(- G \left(\delta
    \left(X_{n-1}\right)\right)\;\1_{\left\{X_{n-1}\in\BB^n\right\}}\Big).
\end{align*}
Hence
$$
\Ex_{\theta
  _0}\left(\left. e^{\left[\frac 1 2 \ell\left(X_{n-1},\varepsilon
 _{n}\right)+G\left(\delta
    \left(X_{n-1}\right)\right)\right]\1_{\left\{X_{n-1}\in\BB^n
    \right\}}}\right|{\cal F}_{n-1}\right)=1
$$
and \eqref{eg} follows.
Now we have
\begin{multline*}
\Pb_{\theta    _0}\left\{\sum_{j=0}^{n-1}\ell\left(X_j,\varepsilon _{j+1}\right)
  \,\1_{\left\{X_j\in \BB^n   \right\}}>-2cu\right\}=\Pb_{\theta    _0}\left\{
  \frac 1 2 S_n^{\left(1\right)}+S_n^{\left(2\right)} - S_n^{\left(2\right)}
  >-cu\right\} \\
\leq \Pb_{\theta    _0}\left\{
  \frac 1 2 S_n^{\left(1\right)}+S_n^{\left(2\right)}
  >\frac 1 2 cu\right\} +\Pb_{\theta    _0}\left\{
   - S_n^{\left(2\right)}
  >-\frac 3 2cu\right\} \leq e^{-\frac 1 2 cu} +\Pb_{\theta    _0}\left\{
    S_n^{\left(2\right)}
  <\frac 3 2 cu\right\}
\end{multline*}
where we used \eqref{eg}. In view of \eqref{une} and \eqref{doux}, it is left to show that for all $p>1$,
$$
\Pb_{\theta    _0}\left\{
    \sum_{j=0}^{n-1}G \left(\delta
    \left(X_{j}\right)\right)\1_{\left\{X_j\in \BB^n
    \right\}}
  <\frac 3 2 cu\right\}\leq \frac{C_p}{\left|u\right|^p}.
$$

Following  \cite{ChK11},
we shall split the consideration into the cases $u<n^s$ and $n^s\leq u<n\left(\beta -\alpha\right)$,
for some $s\in\left(0,1\right)$.

To this end, note that the Hellinger integral $H(\delta)$ is a continuous function of $\delta$:
\begin{align*}
&\big(H(\delta)-H(\delta+\eta)\big)^2  =  \left(\int_{-\infty}^\infty \left[\left(\frac{f(\delta+y)}{f(y)}\right)^{1/2}-
\left(\frac{f(\delta+\eta+y)}{f(y)}\right)^{1/2} \right]f(y) dy\right)^2 \le \\
& \int_{-\infty}^\infty \left( \left(\frac{f(\delta+y)}{f(y)}\right)^{1/2}-
\left(\frac{f(\delta+\eta+y)}{f(y)}\right)^{1/2}\right)^2 f(y) dy =\\
& 2- 2
\int_{-\infty}^\infty   \sqrt{ f(\delta+y) f(\delta+\eta+y)} dy =
\int_{-\infty}^\infty   \Big( \sqrt{f( y)} - \sqrt{f( \eta+y)}\Big)^{2} dy \le \\
& \int_{-\infty}^\infty \big|f( y)-f( \eta+y)\big|dy\xrightarrow{\eta\to 0}0
\end{align*}
where we used LeCam's inequality for the Hellinger and the total variation distances
and the convergence holds by Scheffe's lemma.

By continuity of $G(\delta)=-\ln H(\delta)$  and since $G(\delta)>0$ for all $\delta\ne 0$, the assumption \ref{c2}
implies that for $u<n^s$
\begin{align*}
G\left(\delta
    \left(X_{j}\right)\right)\1_{\left\{X_j\in \BB^n
    \right\}}\geq \inf_{\theta _0\leq v\leq \theta _0+n^{s-1}}G\left(\delta
    \left(v\right)\right)\geq c_0
\end{align*}
with some constant $c_0>0$ and
\begin{align*}
\Pb_{\theta _0}\left\{ S_n^{\left(2\right)} <\frac 3 2cu\right\}\leq
\Pb_{\theta _0}\left\{ \sum_{j=0}^{n-1}\1_{\left\{X_j\in \BB^n \right\}}
  <\frac 3 2\frac{ c}{c_0}\,u\right\}.
\end{align*}
Now let $\eta _j\left(u\right)=\Ex_{\theta
  _0}\1_{\left\{X_j\in \BB^n \right\}}
-\1_{\left\{X_j\in \BB^n \right\}} $. Since the density $f(x)$ is continuous and positive, so is the invariant density $\varphi \left(x,\theta _0\right)$ and
\begin{align*}
  S_n^{\left(3\right)}&=\sum_{j=0}^{n-1}\Ex_{\theta
  _0}\1_{\left\{X_j\in \BB^n\right\}} =n\int_{\theta
   _0}^{\theta _0+u/n}\varphi \left(x,\theta _0\right)\,{\rm d}x\geq C'\,u,
\end{align*}
with a positive constant $C'$.
Then
\begin{multline*}
\Pb_{\theta _0}\left\{ \sum_{j=0}^{n-1}\1_{\left\{X_j\in \BB^n
    \right\}} <\frac 3 2\frac{c}{c_0}\,u\right\} =\Pb_{\theta _0}\left\{
  -\sum_{j=0}^{n-1} \eta _j\left(u\right)
  <-\left(S_n^{\left(3\right)}-\frac 3 2\frac{c}{c_0}u \right) \right\}\\
  \qquad
  \qquad \leq \Pb_{\theta _0}\left\{ \sum_{j=0}^{n-1} \eta _j\left(u\right)
  >\kappa u \right\} \leq \frac{\Ex_{\theta _0} \left|\sum_{j=0}^{n-1} \eta
    _j\left(u\right) \right|^{2p}}{\left|\kappa u\right|^{2p}},
\end{multline*}
where   we chose  $c$ small enough, so that $C'-\frac 3 2c/c_0=\kappa >0$.
Using the geometric mixing property \ref{c4} and an appropriate version of Rosenthal's inequality
as in Lemma 2.2 \cite{ChK11}, we get
$$
\Ex_{\theta _0} \left|\sum_{j=0}^{n-1} \eta
    _j\left(u\right) \right|^{2p}\leq C\left(p\right)\:\left|u\right|^p
$$
which yields \eqref{e2} for $|u|<n^s$.
The complementary case, $n^s\le |u|\le (\beta-\alpha)n$ is treated exactly as in Lemma 2.2, \cite{ChK11} or Lemma 2.5, \cite{CKL1}.

\section{Discussion}

Theorem \ref{T1} can be directly  generalized to the multi-threshold autoregression
$$
X_{j+1}=\sum_{k=0}^{K}h_k\left(X_j\right)\1_{\left\{\theta _{k}<X_j\leq
  \theta    _{k+1}\right\}}+\varepsilon _{j+1},\quad j=0,1,\ldots, n,
$$
where ${\bm \theta} =\left(\theta _1,\ldots,\theta _K\right)$ is the unknown
parameter (and $\theta_0=-\infty$ and $\theta_{K+1}=\infty$ are set).  As in \eqref{nG},  $(\varepsilon _j)$
are independent random variables with
known density  $f\left(x\right)>0,x\in \RR$ and the functions
$h_k\left(\cdot \right) $  and $f\left(\cdot \right)$
are continuous and such that $(X_j)$ is geometrically mixing.
Assume  that $\theta _k\in\left(\alpha _k,\beta _k\right)$, where
$\beta_k<\alpha _{k+1}$.

For all sufficiently large $n$ and $u_k\ge 0$, the normalized likelihood ratio is given by
\begin{align*}
 Z_n\left({\bf
  u}\right)&=\prod_{j=0}^{n-1}\frac{f\left(X_{j+1}-\sum_{k=0}^{K}h_k\left(X_j\right)\1_{\left\{\theta
    _{k}+\frac{u_k}{n}<X_j\leq
    \theta_{k+1}+\frac{u_{k+1}}{n}\right\}}\right)
}{f\left(X_{j+1}-\sum_{k=0}^{K}h_k\left(X_j\right)\1_{\left\{\theta
    _{k}<X_j\leq \theta_{k+1}\right\}}\right)}\\
&=\prod_{j=0}^{n-1}\frac{f\left(\sum_{k=0}^{K}h_k\left(X_j\right)\left[\1_{\left\{\theta
    _{k}<X_j\leq \theta_{k+1}\right\}}
    -
    \1_{\left\{\theta
    _{k}+\frac{u_k}{n}<X_j\leq
    \theta_{k+1}+\frac{u_{k+1}}{n}\right\}}\right]+\varepsilon _{j+1}\right)
}{f\left(\varepsilon _{j+1}\right)}\\
&=\prod_{j=0}^{n-1}\frac{f\left(\sum_{k=1}^{K}\left[h_{k-1}\left(X_j\right)-h_{k}\left(X_j\right)\right]
\1_{\left\{\theta
    _{k}<X_j\leq
    \theta_{k}+\frac{u_{k}}{n}\right\}}+\varepsilon _{j+1}\right)
}{f\left(\varepsilon _{j+1}\right)},
\end{align*}
and thus
\begin{align*}
\ln Z_n\left({\bf u}\right)=&
\sum_{j=0}^{n-1}\ln \frac{f\left(\sum_{k=1}^{K}\left[h_{k-1}\left(X_j\right)-h_{k}\left(X_j\right)\right]
\1_{\left\{\theta
    _{k}<X_j\leq
    \theta_{k}+\frac{u_{k}}{n}\right\}}+\varepsilon _{j+1}\right)
}{f\left(\varepsilon _{j+1}\right)} =\\
&
\sum_{k=1}^{K}
\sum_{j=0}^{n-1}\ln \frac{f\big(\delta_k(X_j)
+\varepsilon _{j+1}\big)
}{f\big(\varepsilon _{j+1}\big)} \1_{\left\{\theta
    _{k}<X_j\leq
    \theta_{k}+\frac{u_{k}}{n}\right\}},
\end{align*}
where $\delta_k(x):= h_{k-1}(x)-h_{k}(x)$.
Using the same approach as in the proof of Theorem \ref{T1}, it can be seen that
\begin{align*}
\ln Z_n\left({\bf u}\right)\quad \Longrightarrow\quad
\sum_{k=1}^{K}\sum_{l=1}^{N_k^+\left(u_k\right)} \ln \frac{f\big(
\varepsilon^+_{k,l}+\delta_k(\theta_k)\big)
}{f\big(\varepsilon^+_{k,l}\big)},
\end{align*}
where $N_k^+\left(u_k\right), u_k\geq 0$ are independent Poisson processes
with intensities $\varphi \left(\theta_k,\theta _k\right)$ and $\eps^+_{k,l}$ are i.i.d. random variables with the density $f$.
Similar asymptotic is obtained for $u_k<0$. Consequently the
limit likelihood ratio is a product on $K$ independent one-dimensional
copies of the process defined \eqref{Z} (with $\theta_0$ replaced by $\theta_k$'s) and the corresponding Bayes estimates $\widetilde \theta_{k,n}$, $k=1,...,K$ are asymptotically
independent with the  asymptotic distribution as in Theorem \ref{T1}.

\section*{Acknowledgement} The authors are grateful to the referee for the careful proofreading of the manuscript and
the suggested improvements.

\end{document}